\documentclass[12pt,reqno,final]{amsart}
\usepackage{amsmath, amssymb, epsfig}
\usepackage{mathrsfs}
\usepackage{color}
\usepackage{alg}
\usepackage{url}
\usepackage{multirow} 

\usepackage[font=small,margin=10pt,labelfont={bf},labelsep={space}]{caption}
\usepackage{subfig}
\usepackage{wrapfig}

\usepackage[scaled]{helvet} 
\usepackage{fourier} 
\usepackage{bm} 

\makeatletter
\def\algspacing{\alg@unmargin}
\makeatother

\newcommand{\BLD}[1]{\mbox{\boldmath $#1$}}

\usepackage[margin=1in]{geometry}

\newlength{\algorithmwidth}
\theoremstyle{plain}

\theoremstyle{definition}

\theoremstyle{remark}

\numberwithin{theorem}{section}
\numberwithin{equation}{section}

\begin{document}

\title[]{On the Mathematics of Music: From Chords to Fourier Analysis}
\author{Nathan Lenssen and Deanna Needell}
\date{\today}

\begin{abstract}
Mathematics is a far reaching discipline and its tools appear in many applications.  In this paper we discuss its role in music and signal processing by revisiting the use of mathematics in algorithms that can extract chord information from recorded music.  We begin with a light introduction to the theory of music and motivate the use of Fourier analysis in audio processing.  We introduce the discrete and continuous Fourier transforms and investigate their use in extracting important information from audio data.
\end{abstract}

\maketitle

\section{Introduction}
Music is a highly structured system with an exciting potential for analysis. The vast majority of Western music is dictated by specific rules for time, beat, rhythm, pitch, and harmony. These rules and the patterns they create entice mathematicians, statisticians, and engineers to develop algorithms that can quickly analyze and describe elements of songs. In this paper we discuss the problem of \textit{chord detection}, where one wishes to identify played chords within a music file.  With the ability to quickly determine the harmonic structure of a song, we can build massive databases which would be prime for statistical analysis.  A human preforming such a task must be highly versed in music theory and will likely take hours to complete the annotation of one song, but an average computer can already perform such a task with reasonable accuracy in a matter of seconds \cite{ellis1}. Here we explore the mathematics underlying such a program and demonstrate how we can use such tools to directly analyze audio files.

\subsection{Organization} In Section~\ref{sec:music} we provide the reader with a brief introduction to music theory and motivate the need for mathematical analysis in chord detection. 
Section~\ref{sec:Fourier} contains an introduction to the mathematics necessary to derive the discrete Fourier transform, which is included in Section~\ref{sec:derive}. 
We conclude with a description of the Fast Fourier Transform and an example of its use in chord detection in Section~\ref{sec:chroma}.

\section{Introduction to Music Theory} \label{sec:music}
We begin with some musical terminology and definitions.

\subsection{Pitches and Scales}
\begin{figure}[ht] 
  \centering
  \includegraphics[width=4in]{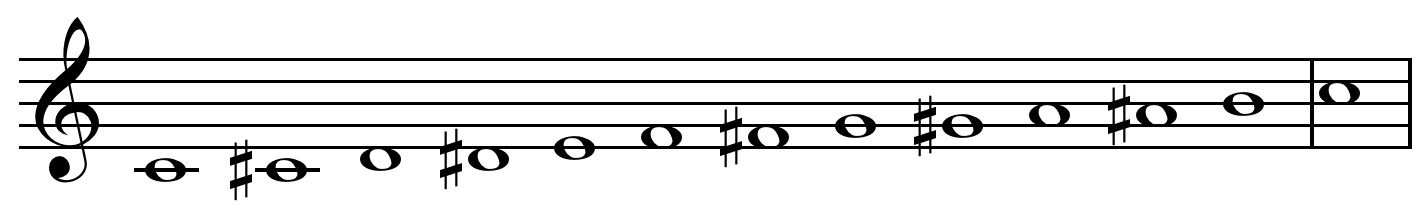}
  \caption{A chromatic scale beginning and ending at C. Notice there are 13 notes because C is played at both the top and bottom.  Image from~\cite{chromaticScale}.\label{fig:chromscalepic}}
\end{figure}

For our purposes we define a \textit{pitch} as the human perception of a sound wave at a specific frequency. For instance, the tuning note for a symphony orchestra is A4 which has a standardized frequency of 440Hz. In the notation A4, A indicates the chroma or quality of the note while 4 describes the octave or height. 
For ease of notation, a \textit{scale} is a sequence of pitches with a specific spacing in frequency. As we follow the pitches of a scale from bottom to top, we start and end on the same note one octave apart (e.g. from C3 to C4). Pitches an octave apart sound similar to the human ear because a one octave increase corresponds to a doubling in the frequency of the sound wave. Western music uses the chromatic scale in which each of the 12 chroma are ordered over an octave. These 12 notes are spaced almost perfectly logarithmically over the octave. We can use a recursive sequence~\cite{kLee} to describe the chromatic scale: 
\begin{equation}
P_i=2^{1/12}P_{i-1},
\end{equation}
where $P_i$ denotes the frequency of one pitch, and $P_{i-1}$ the frequency of the previous. We can hear the chromatic scale by striking every white and black key of a piano in order up an octave or visualize it by a scale such as that depicted in Figure \ref{fig:chromscalepic}.  Note that pitches use the names A through F, along with sharp ($\sharp$) and flat ($\flat$) symbols, see~\cite{completeMusician} for more on musical symbols.

By using the chromatic scale as a tool, we can  construct every other scale in Western music through the use of intervals. An \textit{interval} refers to a change in position along the 12 notes of the chromatic scale.  We define the interval of a ``half step" or H as a change in one pitch along the chromatic scale. Two half step intervals makes a ``whole step" and is denoted by W. We also use interval of three half steps known as an ``augmented second" denoted by A. Scales are defined by a sequence of these intervals with the condition that the total sum of  steps must equal 12. This guarantees that we start and end on the same chroma known as the root R. There are four prevalent scales in Western music: major, minor, augmented, and diminished. The intervals that describe these scales can be found in Table 2.1. The table is used by selecting any starting note as the root and then using the intervals to construct the remaining notes. For instance, a C minor (Cm) scale is C,D,E$\flat$,F,G,A,B$\flat$,C. 
\begin{table}[ht] 
\begin{center} 
    \begin{tabular}{ |l|l| }
    \hline
    {\bf Scale ``Color"} &  {\bf Defining Steps} \\ \hline
    Chromatic & RHHHHHHHHHHHHR \\ \hline
    Major & RWWHWWWHR \\ \hline
   (Natural) Minor & RWHWWHWWR \\ \hline
    Diminished & RHWHWHWHWR \\ \hline
    Augmented &  RAHAHAHR\\ 
    \hline
    \end{tabular}
\caption{Interval construction of the four core scales with the chromatic scale for reference. Note that the intervals apply for when ascending in pitch only. When determining the descending scale, the order of intervals is reversed.}  
\end{center}
\end{table}

\subsection{Triads and Chords}
Multiple pitches played simultaneously are defined as a chord. Chords are essential in music analysis because they compactly describe the entire melodic and harmonic structure of a section of music. That is, a chord indicates what notes and combination of notes should be played at a moment in time. A \textit{triad} is a specific and simple chord containing the first, third and fifth note of a scale (I, III, V).
\begin{figure}[t] 
  \centering
  \includegraphics[width=3.5in]{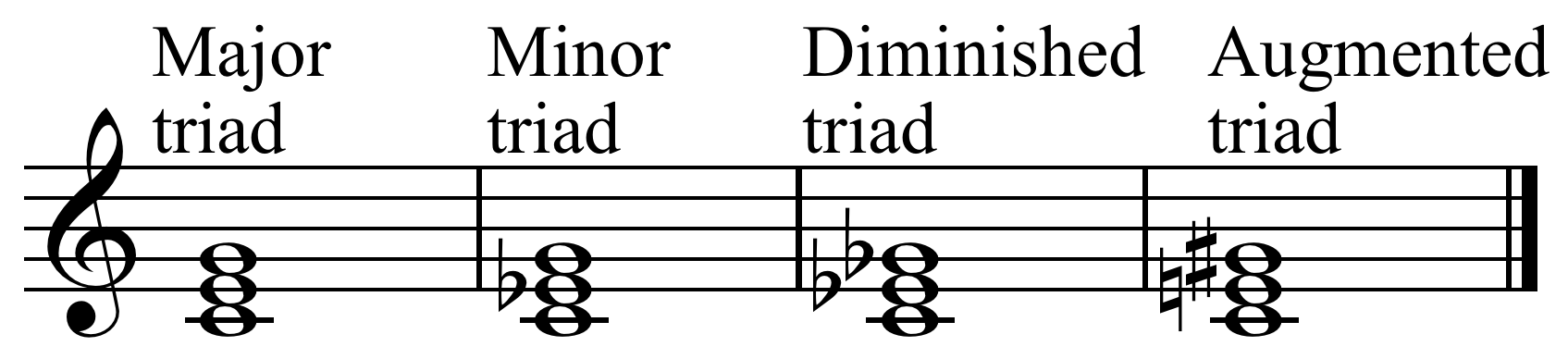}
  \caption{Triads in the key of C. Image from \cite{triad}.} \label{fig:ctriads}
\end{figure}
Figure \ref{fig:ctriads} shows the triads for each of the four scales in the key of C. These triads describe the major (maj), minor (min), diminished (dim), and augmented (aug) chords in our model. 

While triads are a useful way to understand the tonal structure of music, four notes are often needed to completely describe tonal character. Adding and possibly altering the seventh note of the scale (VII) creates new and essential chords. The three chords that we need a seventh to describe are the major seventh (maj7), minor seventh (min7), and dominant seventh (dom7) chords. The major and minor seventh chords follow directly from the major and minor scales. They each contain the I, III, V, VII of their respective scales. The dominant seventh chord does not follow one of the scales we have described. With respect to the major scale it contains the I, III, V, VII$\flat$. The theory behind the dominant seventh chord is a consequence of the theory of musical modes, which we refer the interested reader to~\cite{completeMusician} for more information.


\subsection{Chord Inversions}
We have described 7 chord families and 21 roots giving a total of 147 different possible chords.  However, this number is assuming that the root is always the lowest note in the chord, which isn't always the case.  Instead, a \textit{chord inversion} may be used. An inversion can be described algorithmically as follows: first, the root is raised an octave. This is known as the first inversion. If we repeat that process again (with the new lowest note), we are left with the second inversion. For a triad, two inversions recover the original chroma arrangement. When chords involve a seventh, one can perform three inversions. Including inversions, the total number of chords is
\begin{equation}
C = (21 \text{ roots})[(4 \text{ triads})(2 \text{ inversions}) + (3\;\; 7^{th}\,\text{chords})(3 \text{ inversions})] = 357.
\end{equation}
Distinguishing between this amount of distinct possible chords is quite a task, particularly since the octave information of a note cannot help us narrow down the chord choice. In the next section we will explore how mathematicians and engineers take raw audio files and determine the chord progressions.
\begin{figure}[t] 
  \centering
  \includegraphics[width=2.3in]{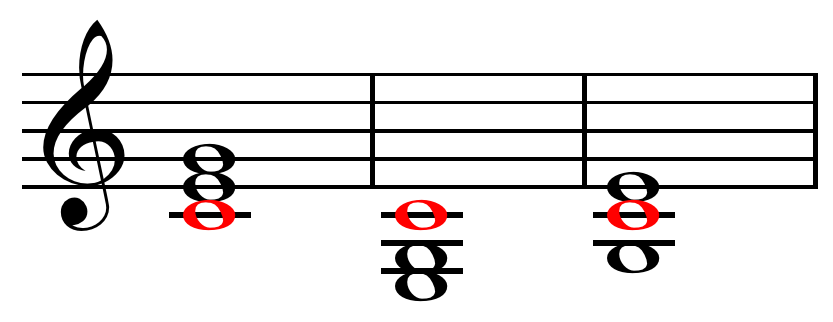}
  \caption{The root and inversion triads in the key of C\cite{inversion}.} \label{fig:cinversiontriads}
\end{figure}

\subsection{Towards Chord Detection}
Although chord analysis has been widely studied (see e.g.~\cite{bartsch2001catch,fujishima1999realtime,gold2011speech,kLee,raphael2002automatic,ellis1} and the references therein), we will focus on a chord recognition model developed by Sheh and Ellis in 2003~\cite{ellis1}. This model can be broken down into two major components: signal analysis and chord fitting. The goal of the signal analysis portion of the model is to break down a raw audio signal into chroma intensities. These intensities are then fed into the chord fitting portion where the most likely chord representations for the music are determined.

The primary mathematical tool used to decompose the raw audio signal is the Fast Fourier Transform (FFT), which is an optimized discrete Fourier transform algorithm. The FFT is used to determine the fundamental frequencies and therefore pitches that are present in the raw signal. However, the octave of the pitch is generally irrelevant to the chord identity, so one needs to transform the pitches obtained through the FFT into octave-independent chroma. The chroma information of an audio file is known as the Pitch Class Profile (PCP)~\cite{fujishima1999realtime}.

Once the PCP of the input is discovered, it is used to determine the best fitting chords for the music. To do this, a Hidden Markov Model that has been trained by a large set of pre-annotated audio files is used. The optimal chord assignments given a PCP are determined through the expectation maximization (EM) algorithm\cite{gold}. 
The Viterbi alignment algorithm~\cite{gold} can be used to forcibly align the chord labels with the timing of the music, creating a fully analyzed track. 

We now explore the FFT by building an understanding of the underlying mathematics. In doing so, we reveal how Fourier analysis is useful in determining the chroma structure of an audio file.

\section{Introduction to the Fourier Transform}\label{sec:Fourier}

\subsection{Frequency and Time Domains} \label{sec:Fintro}
The development of the mathematical construct now referred to as the \textit{Fourier transform} was motivated by two problems in physics that are very prevalent and observable. These problems are the conduction of heat in solids and the motion of a plucked string whose ends are fixed in place \cite{skinnybook,briggs}. One instantly sees the relevance to chord detection since string instruments such as violins or pianos produce sound by amplifying the vibrations of a fixed string. The history of the development of the Fourier series and transform is interesting and rich. A reader interested in its development is recommended to read the introductions of \cite{skinnybook,briggs}. History aside, the most crucial information from the development of Fourier analysis is that functions can be represented in both the \textsl{time domain} as well as the \textsl{frequency domain}. 

The connection between the two domains is easiest to see in a periodic system such as a vibrating string. One would likely describe the motion of a string by focusing on the changing position of points on the string over time. More rigorously, one can model such motion using a differential equation where the initial condition is the initial displacement~\cite{skinnybook}. Solving such differential equations yields a function $f(t)$ (where $t$ is time) which represents the motion of the string. This function $f(t)$ is known as the \textit{time domain representation} of the motion of the string; it provides information about how the string behaves as time progresses.

\begin{figure}[t] 
  \centering
  \includegraphics[width=2.5in]{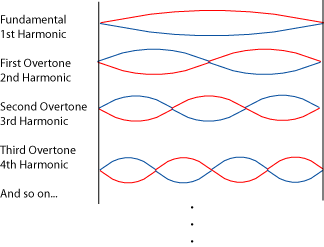}
  \caption{The first four standing waves of a vibrating, fixed string. Image from \cite{swave}.} \label{fig:swave}
\end{figure}

Now imagine one plucks the string in precisely a way that makes the string vibrate only at the fundamental harmonic illustrated in Figure \ref{fig:swave}.  This pattern is represented in the time domain by a single sinusoid with frequency $\nu_0$.  Notice that the motion of the string is completely described by this frequency $\nu_0$ and the amplitude of the oscillation.
Thus, the frequency domain $F(\nu)$ of this specific string has just one spike at $\nu=\nu_0$ with the height of the spike equal to the amplitude of the wave. The example of the fundamental harmonic helps to visualize the frequency domain, but in real systems there is typically more than one frequency. To account for this, one constructs the frequency domain representation by an infinite series of these harmonics weighted in such a way that they represent the motion of the string. This series is known as the \textit{Fourier series} and provides the basis for the Fourier transform. Through the Fourier transform one is able to obtain the frequency domain representation of a time domain function. The Fourier transform is invertible with the inverse Fourier transform returning the time domain function from a frequency domain function.

Audio signals are recorded and listened to in the time domain; they contain the intensity of the sound as a function of time. However, pitches and therefore chords are represented by frequencies as shown in Section~\ref{sec:music}. Therefore, the Fourier transform is used to convert the time domain input signal into a frequency representation which can be analyzed for intensities of specific pitches. The relationship between the pitches at a given time is then analyzed to determine the chord that is being played. It is important to note that the Fourier transform is a transformation in the complex plane; we will see that even for real-valued inputs such as audio signals, the Fourier transform will return a complex-valued frequency domain representation.

\subsection{The Continuous Fourier Transform}
We first define $\omega_k$ as the angular frequency through its relationship with the ordinary frequency $\nu$ by
\begin{equation}
\omega_k\equiv 2\pi k \nu
\end{equation} 

The relationship between the time domain function $f(t)$ and its corresponding frequency domain function $F$ is defined by the \textit{Fourier transformation}:
\begin{equation} \label{eq:fourier}
F(\omega_k) \equiv \int^{\infty}_{-\infty}\!f(t)\,e^{-2 \pi i k t} \, \mathrm{d} t \qquad k\in (-\infty,\infty)
\end{equation}
The presence of the complex sinusoidal is clear by recalling that $e^{i\omega t} = \cos (\omega\, t) + i\sin (\omega\, t)$. From~\eqref{eq:fourier} one can derive the inverse Fourier transform:
\begin{equation}
f(t) = \int^{\infty}_{-\infty}\!F(\omega_k)\,e^{2 \pi i k t} \, \mathrm{d} k \qquad k\in (-\infty,\infty)
\end{equation}

Fourier analysis in the continuous case is a rich subject of study.  However, digital applications require its discrete counterpart, the discrete Fourier transform (DFT).


\subsection{The Discrete Fourier Transform}
Motivated by the continuous case, one similarly defines the discrete Fourier transform for a complex vector $f$ centered at 0 by:

\begin{equation} \label{eq:DFT}
F_k \equiv \sum_{n=-\frac{N}{2}+1}^{\frac{N}{2}} f_n e^{-i2 \pi nk/N} \qquad k=0:N-1,
\end{equation}
where we write $f_n$ to denote the $n$th entry of the vector $f$.
Likewise, the inverse DFT has the form
\begin{equation}
f_k = \sum_{n=-\frac{N}{2}+1}^{\frac{N}{2}} F_n e^{i2 \pi nk/N} \qquad k=0:N-1.
\end{equation}
We next study these transforms in more detail, with the focus of chord detection and audio processing.

\subsection{Sinusoids}
We define a sinusoid as a function of the form
\begin{equation}
x(t) = A\sin(2\pi \nu t + \phi)
\end{equation}
When discussing audio signals, as in~\cite{stanford} we let
\begin{align*}
A &= \mbox{Amplitude} \\
2\pi \nu &= \mbox{frequency(Hz)}\\
t &= \mbox{time (s)} \\
\phi &= \mbox{initial phase (radians) }\\
2\pi \nu t + \phi &= \mbox{instantaneous phase (radians). } 
\end{align*}
Fourier transforms are built on the complex properties of sinusoids which follow from Euler's identities,
\begin{align}
e^{i\theta} &= \cos(\theta) + i\sin(\theta) \label{eq:euler1} \\
e^{\pm i 2 \pi \nu x} &= \cos(2 \pi \nu x) \pm i \sin(2 \pi \nu x) \label{eq:euler2}, 
\end{align}
the latter being of the form most relevant to audio.

A beautiful and useful fact is that sinusoids are orthogonal (for a proof of this, see for example Section 3 of~\cite{nathan}).  This allows us to utilize them as an orthonormal basis to represent vectors like audio signals.

\subsection{The Delta ($\delta$) Function} \label{sec:deltafunction}
Although sound is naturally analog, computers must store sound in digital format.  For that reason, sound must be \textit{sampled} in an efficient way for an accurate representation.
The mathematics of sampling is a large field and we refer the reader to a work such as~\cite{pharr2010physically} for a thorough discussion, and~\cite{CSwebpage,griffin2011single} for new work in this area. Given a continuous signal $f(t)$ as an input over time $t\in[0,T]$, the sampler returns a time-series vector that contains values according to some sampling frequency with a uniform gap $\Delta t$ between samples. This vector can be viewed as a continuous ``function" that is a piece-wise collection of Dirac delta functions $\delta$, each centered at some multiple of our sampling period $n\Delta t$. We say ``function" because $\delta$ is not a function by the rigorous definition, but rather a distribution or generalized function~\cite{briggs}. The delta function, or \textit{spike}, is a special mathematical construct with a number of properties that are intuitive and useful to the study of the DFT. We now outline a few of its most relevant properties to this topic.

First, the $\delta$ function is non-zero at every point but 0 at which it is infinite. That is,
\begin{equation}
\delta(t)=
\begin{cases}
\infty & \text{if } t=0 \\
0 & \text{if } t \neq 0.
\end{cases}
\end{equation}
This property illuminates why we use the term spike. The $\delta$ function is an infinitely thin, infinitely tall peak centered around zero.  A $\delta$ function centered at $t=a$ can of course be written as $\delta(t-a)$. 

Second, the area under the delta function is defined to be 1. That is,
\begin{equation}
\int_{-\infty}^{\infty} \delta(t) \, \hbox{d}t=1
\end{equation}
This property can be exploited to reveal the useful general sifting property to discretize integrals,
\begin{equation}
\int_{-\infty}^{\infty} f(t) \delta(t-a) \, \hbox{d}t=f(a)
\end{equation}
The relevance of the general sifting property to sound processing cannot be overstated. 
For a continuous input signal $f(t)$, the sifting property and $\delta(t-a)$ can be used to obtain an evenly spaced sequence of $a$ values corresponding to the instantaneous values of the original function. While certainly not a sophisticated method, it provides a relevant and simple application that demonstrates the power of the $\delta$ function. 


\subsection{Fourier Transforms of $\delta$ Functions}
We begin by looking at the Fourier transform of a single spike centered at zero, which will be instrumental in our derivation of the DFT. By the general sifting property,
\begin{equation}
F(t) = \int_{-\infty}^{\infty} \delta(t)\, e^{-i2\pi \nu t}\, \hbox{d}t \medspace = \medspace e^0\medspace = \medspace 1.
\end{equation}
This means that a spike in the time domain translates to a flat line in the frequency domain (an even weight among all frequencies). This is a somewhat physically curious result  that can be reconciled by understanding that for every frequency, $\cos (2\pi \nu \cdot0)=1$ and at every other value of $t$, the sinusoids take on every value between 0 and 1. 

We also see an interesting result when we have a spike at the origin in the frequency domain. Intuitively this spike symbolizes that the time dependent function is made up of a single sinusoid with frequency zero, or a straight line. We can confirm this by looking at the inverse Fourier transform of a $\delta$ function
\begin{equation}
f(t) = \int_{-\infty}^{\infty} \delta(t)\, e^{i2\pi \nu t}\, \hbox{d}t \medspace = \medspace e^0\medspace = \medspace 1.
\end{equation}
These results are useful and interesting, but do not provide us with the connection between a spike and a sinusoid necessary to build the DFT. It appears that a spike in the frequency domain at any location besides the origin will correspond to some sort of sinusoid in the time domain. We can see this is true by taking the inverse Fourier transform of a $\delta$ function shifted by some frequency $\nu_0$ which yields
$$
f(t) = \int_{-\infty}^{\infty} \delta(\nu-\nu_0)e^{i2\pi \nu x} \, \hbox{d}x  =  \medspace e^{i2\pi \nu_0 x}  
      =  \cos (2\pi \nu_0 x) + i \sin(2\pi \nu_0 x).
$$
By clever addition of pairs of Inverse Fourier transforms, we can determine what spikes in the frequency domain have inverse Fourier transforms that are real-valued sine and cosine functions.  We notice that to get a sine function in the time domain we must construct spikes in such a way that their sine components cancel. Recalling the definition of cosine from Euler's formula we have,
\begin{align}
f(t)  \medspace &=  \medspace \int_{-\infty}^{\infty} \left [ \frac{1}{2} \delta(\nu-\nu_0) + \frac{1}{2} \delta(\nu+\nu_0) \right ] e^{i2\pi \nu x} \, \hbox{d}x \nonumber \\
&= \medspace  \frac{1}{2}\int_{-\infty}^{\infty} \delta(\nu-\nu_0)e^{i2\pi \nu x} \, \hbox{d}x + \frac{1}{2}\int_{-\infty}^{\infty} \delta(\nu+\nu_0)e^{i2\pi \nu x} \, \hbox{d}x \nonumber \\
&= \medspace  \frac{e^{i2\pi \nu_0 x} + e^{-i2\pi \nu_0 x}}{2} \nonumber \\
&= \medspace  \cos(2\pi \nu_0 x)
\end{align}
Similarly for the sine function we have the inverse Fourier transform
\begin{align}
f(t)  \medspace &=  \medspace \int_{-\infty}^{\infty} \left [ \frac{1}{2i} \delta(\nu-\nu_0) - \frac{1}{2i} \delta(\nu+\nu_0) \right ] e^{i2\pi \nu x} \, \hbox{d}x \nonumber \\
&= \medspace  \frac{1}{2i}\int_{-\infty}^{\infty} \delta(f-f_0)e^{i2\pi \nu x} \, \hbox{d}x - \frac{1}{2i}\int_{-\infty}^{\infty} \delta(\nu+\nu_0)e^{i2\pi \nu x} \, \hbox{d}x \nonumber \\
&= \medspace  \frac{e^{i2\pi \nu_0 x} - e^{-i2\pi \nu_0 x}}{2i} \nonumber \\
&= \medspace  \sin(2\pi \nu_0 x).
\end{align}
In other words, the cosine function of frequency $\nu$ has a Fourier transform of two positive real-valued spikes at $\pm \nu$ in the frequency domain. A sine function of frequency $\nu$ has a Fourier transform that lies purely in the imaginary frequency domain with a negative spike at $+\nu$ and a positive spike at $-\nu$. The sine function serves as a demonstration of the necessity of complex numbers in the Fourier transform; a real valued sine wave is described by a completely imaginary frequency representation.

\section{A Derivation of the Discrete Fourier Transform}\label{sec:derive}
\subsection{Spike Trains and the Discrete Fourier Transform} \label{sec:devDFT1}
As we have defined a spike as a $\delta$ function, a \textit{spike train} is simply a collection of $\delta$ functions. 
For application to the DFT, we fix the spacing of these spikes in time as $\Delta t$. Thus the $n^{\text{th}}$ spike is located at $t_n = n\Delta t$. We define the spike train over time $h(t)$ as
\begin{equation}
h(t) = \sum_n f_n \delta (t-t_n).
\end{equation}
In this form the vector $f$ represents the intensity or magnitude of the spikes. By recalling the general sifting property of the delta function, we can rewrite a spike train as,
$$
h(t) = \sum_n f_n \delta (t-t_n) = \sum_n f(t_n) \delta (t-t_n)  
					       = \sum_n f(t) \delta (t-t_n).
$$
Thus, we can easily construct the spike train of a function by summing a series of evenly spaced delta functions multiplied by the original continuous function. 

With the definition of a spike train we can begin to reveal the DFT. 
Assume we have an evenly spaced spike train of length $N$ on an time interval of length $L$. We center the function at zero so that the interval covered is $[-L/2,L/2]$. Since the points are evenly space in time, we know that the $n^{th}$ point is located at $t_n=nL/N$ for $n=-N/2+1\ldots N/2$. 
Again making use of the general sifting property of $\delta$ functions we have,
\begin{align} \label{eq:almostDFT}
F(\nu) &= \int_{-\infty}^{\infty} \left ( \sum_{n=-\frac{N}{2}+1}^{\frac{N}{2}} f_t \delta (t-t_n) \right ) e^{-i2 \pi \nu t} \, \hbox{d} t \nonumber \\
       &= \sum_{n=-\frac{N}{2}+1}^{\frac{N}{2}} f_n \int_{-\infty}^{\infty} \delta (t - t_n) e^{-i2 \pi \nu t} \, \hbox{d} t \nonumber \\
       &= \sum_{n=-\frac{N}{2}+1}^{\frac{N}{2}} f_n  e^{-i2 \pi \nu x_n}.
\end{align}
Thus the Fourier transform of a spike train is simply the sum of exponentials weighted by their intensities in the spike train \cite{briggs}. We now have an expression that is very close to the definition of the DFT given in~\eqref{eq:DFT}. To make~\eqref{eq:almostDFT} match, we need to define and use the reciprocity relations of the discrete Fourier transform.

\subsection{Reciprocity Relations of the Discrete Fourier Transform}
To uncover the reciprocity relations we will use the same sequence in the time domain as used in the previous section: a vector of length $N$ evenly spaced over a total time of $L$, centered at $t=0$ so that the time interval runs over $[-L/2,L/2]$. If the temporal spacing is $\Delta t$ then the points are defined as $t_n=n\Delta t$. We will use a similar method in discretizing the frequency domain. We again use a vector of length $N$ evenly spaced points over length $\Omega$, again centered at $\nu=0$ so that the interval of frequencies is $[-\Omega / 2, \Omega /2]$. If the spacing of frequencies is $\Delta \nu$ then the points are defined as $\nu_k = k \Delta \nu$. We search for the reciprocity relations that relate the time and frequency domains using the parameters $\Delta t$ and $L$ for the time domain and $\Delta \nu$ and $\Omega$ for the frequency domain. 

Since we only have a discrete number of points in the frequency domain, we are limited in the number of sinusoids we can use to represent the time domain function. Likewise, we are looking at the time domain function over a finite and well-defined time interval $L$. Clearly the longest sinusoid we can resolve in this amount of time is one with $L$ as its period or the fundamental harmonic. Frequency is the reciprocal of the period and we will call this longest period the \textit{fundamental frequency}, which is also the step size in the frequency domain. That is,
\begin{equation} \label{eq:flrelation}
\Delta \nu = \frac{1}{L}.
\end{equation}
Thus the frequencies we will recognize will all be multiples of $\Delta \nu$ and have integer periods over the time interval. It is then clear that the length described in the frequency domain is just the number of points multiplied by the frequency step, or $\Omega = N \Delta \nu$. Using this in conjunction with  \eqref{eq:flrelation} we have the first reciprocity relation,
\begin{equation} \label{eq:rr1}
\Omega = N \Delta \nu = \frac{N}{L} \quad \Longrightarrow \quad {A\Omega = L}
\end{equation}
This relation shows that the lengths of the temporal and frequency domains are inversely proportional and jointly fixed with respect to the input vector length $N$. In practice, we interpret \eqref{eq:rr1} by noting that taking temporal data over a longer range of time means that the DFT yields a smaller range of frequencies and vice versa \cite{briggs}. 

Recalling from \eqref{eq:flrelation} that $\Delta \nu = 1/L$, we can present the second reciprocity relation as
\begin{equation}
\frac{1}{\Delta \nu} = L = N \Delta t \quad \Longrightarrow \quad {\Delta t \Delta \nu  = \frac{1}{N}}.
\end{equation}
This relation reveals very similar information to the first. The spacings of the temporal and frequency domain are inversely proportional and fixed by the number of points or the length of the vector. 

We now have the tools necessary to derive the precise form of the DFT.


\subsection{The Discrete Fourier Transform}
Recall from \eqref{eq:almostDFT} that
\begin{equation} \label{eq:almostalmostDFT}
F(\nu) = \sum_{n=-\frac{N}{2}+1}^{\frac{N}{2}} f_n  e^{-i2 \pi \nu x_n}.
\end{equation}
Using the assumption that we have $N$ points over a time length of $L$ we know that the step spacing in the frequency domain must be $\nu_k = k \Delta \nu = k/L$ for $k = -L/2+1\ldots L/2$. We use this to show that \eqref{eq:almostalmostDFT} can be rewritten as
\begin{equation} \label{eq:finalDFT}
F(\nu) = \sum_{n=-\frac{N}{2}+1}^{\frac{N}{2}} f_n  e^{-i2 \pi \nu_k x_n} 
            = {\sum_{n=-\frac{N}{2}+1}^{\frac{N}{2}} f_n  e^{-i2 \pi n k /N}}
\end{equation}
The DFT as defined in~\eqref{eq:DFT} is now apparent. 
We proceed by investigating how we employ the DFT and apply it to music signal analysis.

\section{The Chromagram}\label{sec:chroma}
Since the DFT is a sum of indexed values in the exponential kernel, we can express it as the linear equation
\begin{equation}
\BLD{F} = \BLD{W} \BLD{f},
\end{equation}
where $\BLD f$ is the input data in the time domain, $\BLD F$ is the output in the frequency domain, and $\BLD W$ is the nonsingular matrix
\begin{equation}
\BLD{W} =
 \begin{pmatrix}
  e^{-i2 \pi (0)/N} & e^{-i2 \pi (0)/N} & e^{-i2 \pi (0)/N} & \cdots & e^{-i2 \pi (0)/N} \\
  e^{-i2 \pi (0)/N} & e^{-i2 \pi (1)/N} & e^{-i2 \pi (2)/N} &\cdots & e^{-i2 \pi (N-1)/N} \\
  e^{-i2 \pi (0)/N} & e^{-i2 \pi (2)/N} & e^{-i2 \pi (4)/N} &\cdots & e^{-i2 \pi (2(N-1))/N} \\
  \vdots  & \vdots  &  \vdots & \ddots & \vdots  \\
  e^{-i2 \pi (0)/N}  & e^{-i2 \pi (N-1)/N} & e^{-i2 \pi (2(N-1))/N} &\cdots & e^{-i2 \pi ((N-1)(N-1))/N}
 \end{pmatrix}.
\end{equation}
Since $\BLD W$ is nonsingular, we may express the inverse DFT by
\begin{equation}
\BLD{f} = \BLD{W} ^{-1}\BLD{F}
\end{equation}

Through using this method of matrix multiplication, we can calculate the DFT in $O(n^2)$ time where $n$ is the length of the input vector. However, with sampled music, the inputs are extremely high dimensional and we would like to find a method that computes the DFT in a faster amount of time.

\subsection{The Fast Fourier Transform}
\begin{figure}[t]
  \centering
  \includegraphics[width=6in]{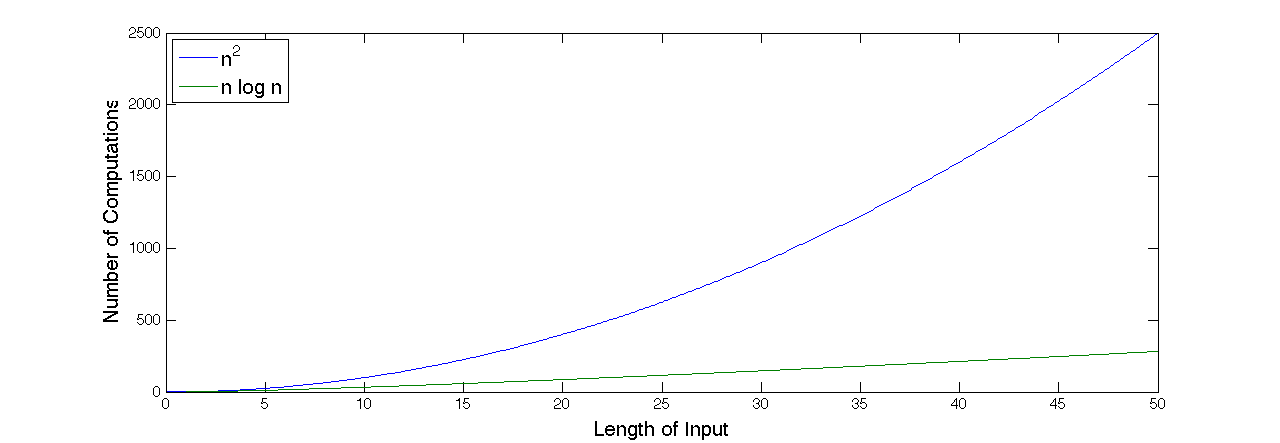}
  \caption{A quick comparison of $O(n^2)$ and $O(n \log_2 n)$ processing speed.} \label{fig:order}
\end{figure}

In 1965, Cooley and Tukey published an algorithm that fundamentally changed the digital signal processing landscape~\cite{cooley}. By exploiting symmetries of the DFT, they were able to reduce the running time of DFTs from $O(n^2)$ to $O(n \log n)$. This algorithm is the first fast Fourier transform (FFT), named for this increase in computational speed. As can been seen in Figure \ref{fig:order} this reduction in processing time is quite significant even for an input vector of length 50. In our examples, we use audio that has sampling frequencies of 11025 Hz. Thus in a three minute song, there are about 2,000,000 input points. In this case, the $O(n \log n)$ FFT algorithm provides a frequency representation of our data 
\begin{equation}
\frac{n^2}{n \log_2 n} = \frac{\left ( 2 \cdot 10^6 \right )^2}{\left ( 2 \cdot 10^6 \right ) \log_2 \left ( 2 \cdot 10^6 \right )} \approx  100,000 \text{ times faster.} \nonumber
\end{equation}
Clearly this boost in speed that the FFT provides is substantial. We can calculate a FFT in MATLAB in fractions of a second when a full DFT would take hours. Thus the FFT facilitates the spectral analysis or frequency domain analysis of large data such as audio files. 

\subsection{Spectrograms and Chromagrams}
By looking at the FFT of overlapping segments of time in an audio file, we can construct the Short-Time Fourier Transform (STFT)~\cite{jacobsen2003sliding}. Each of the FFTs in the STFT reveal the frequency domain representation or spectrum of a small time interval of the input signal. When combining all of these time segments with the proper synthesis, we can construct data on the intensities of frequencies as time progresses \cite{stanford}. The STFT allows us to add a time dimension to the DFT which enables us to observe how the frequency domain changes over time. Since frequencies are representative of pitches, we can use the spectrograms to determine the chroma played at a moment of time in the signal. The creation of these spectrographs is crucial in determining how chords are changing in music. As stated previously, the octave information of the sound is irrelevant. We determine the chroma intensity by collecting all intensities of a note regardless of its octave. The algorithm of collecting these chroma is referred to as a \textit{chromagram}, and it contains the information that we pass to the Hidden Markov Model to determine the chord represented~\cite{ellis1}. 
\begin{figure}[t]
  \centering
  \includegraphics[width=6in]{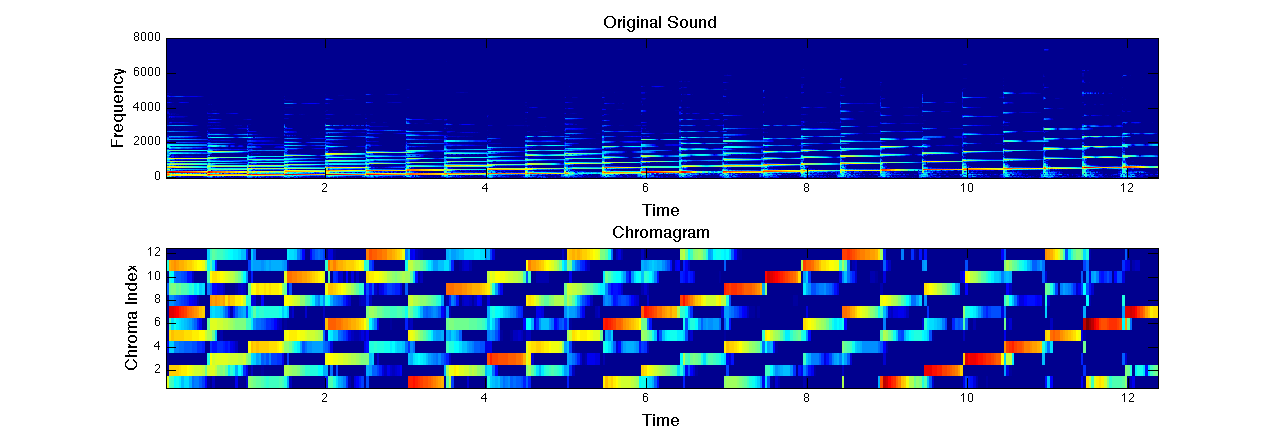}
  \caption{The spectrogram (top) and chromagram (bottom) of an ascending scale.} \label{fig:spec}
\end{figure}\noindent

In this example we use MATLAB code by Ellis \cite{ellis2}. With this code, we are able to create a spectrograph of a $\sim 12$ second audio clip of a piano playing an ascending chromatic scale. The results of our analysis are demonstrated in Figure \ref{fig:spec}. The spectrogram of the sound is depicted on top, from which it is clear when each note is struck as well as the general upward trend in frequency as the scale ascends. Through an algorithm developed by Ellis we are able to extract the chromagram from the spectrogram \cite{ellis2}. 

In the bottom plot of Figure \ref{fig:spec} we can see when each key is struck as well as the upward trend. However, we observe a behavior that when reaching the 12$^\text{th}$ chroma, the major intensity block jumps back down to the first. This jump is clear around $t=9$ and is explained by our disregard for octave information when discussing chroma. Since we are only describing the note name, a jump from B to C is a jump from 12 to 1 in the chroma space. Clearly, the chromagram is accurately describing the pitch characteristics of a chromatic scale and we have the first step in our chord detection model.

\section{Conclusion}
In this paper we have laid a foundation for understanding the mathematics behind a chord recognition model. We have provided the reader with general knowledge of music and a description of the current methods for chord recognition. In addition, we derived the DFT from the Fourier transform and gained an appreciation for the applications of Fourier analysis in signal processing. Finally, we have shown through examples how the FFT of a function can be used to create a chromagram from an audio file.    

\section*{Acknowledgments}

We would like to thank Professor Dan Ellis for his very helpful guidance and suggestions.  In addition, we thank Professor Mark Huber for directing this discussion to the venue of Humanistic Mathematics.

\bibliographystyle{plain}
\bibliography{MathThesisBib}


\end{document}